%% file: incexpm_arxiv.tex
\documentclass[a4paper]{article}
\usepackage[margin=4cm]{geometry}

\usepackage{amssymb,amsmath,amsfonts}
\usepackage{amsthm} 
\usepackage{graphicx} 
\usepackage{lmodern}
\usepackage{url}
\usepackage{pdfsync}
\usepackage{paralist}

\usepackage{algorithm}
\usepackage{algorithmic}

\newcommand{\R}{{\mathbb R}}
\newcommand{\N}{{\mathbb N}}

\newcommand{\bigO}{{\mathcal O}}
\makeatletter
\def\BState{\State\hskip-\ALG@thistlm}
\makeatother
\DeclareMathOperator{\tr}{Tr}

\newcommand{\norm}[1]{\lVert #1 \rVert}

\newcommand{\defby}{\mathrel{\mathop:}=}
\newcommand{\bydef}{=\mathrel{\mathop:}}

\newtheorem{proposition}{Proposition}[section]
\newtheorem{lemma}[proposition]{Lemma}

\newtheorem{remark}[proposition]{Remark}

\title{Incremental computation of block triangular matrix exponentials with application to option pricing}

\author{Daniel Kressner\footnotemark[1]  \and
Robert Luce\footnotemark[1]  \and
Francesco Statti\footnotemark[1]{\hspace*{4pt}}\footnotemark[2]}

\begin{document}

\maketitle

\renewcommand{\thefootnote}{\fnsymbol{footnote}}

\footnotetext[1]{\'Ecole Polytechnique F\'ed\'erale de Lausanne,
Station 8, 1015 Lausanne, Switzerland~({\tt
\{daniel.kressner, robert.luce, francesco.statti\}@epfl.ch})}

\footnotetext[2]{Research supported through the European Research
Council under the European Union's Seventh Framework Programme
(FP/2007-2013) / ERC Grant Agreement n.~307465-POLYTE.}

\begin{abstract}
    \input{abstract}
\end{abstract}


\input{chapter1}
\input{chapter2}
\input{chapter3}

\input{chapter4}
\input{chapter5}

\bibliographystyle{plain}
\bibliography{references}

\end{document}

%% file: abstract.tex
We study the problem of computing the matrix exponential of a block
triangular matrix in a peculiar way: Block column by block column,
from left to right.  The need for such an evaluation scheme arises
naturally in the context of option pricing in polynomial diffusion
models.  In this setting a discretization process produces a sequence
of nested block triangular matrices, and their exponentials are to be
computed at each stage, until a dynamically evaluated criterion allows
to stop.  Our algorithm is based on scaling and squaring.  By
carefully reusing certain intermediate quantities from one step to the
next, we can efficiently compute such
a sequence of matrix exponentials.

%% file: chapter1.tex
\section{Introduction}


We study the problem of computing the matrix exponential for a sequence of nested block
triangular matrices. In order to give a precise problem formulation, consider a sequence of
block upper triangular matrices $G_0, G_1, G_2, \dotsc$ of the form
\begin{equation}
    \label{eq:G_n}
    G_n =
    \begin{bmatrix}
        G_{0,0} & G_{0,1}  & \cdots & G_{0,n} \\
                & G_{1,1}  & \cdots & G_{1,n} \\
                &          & \ddots & \vdots  \\
                &          &        & G_{n,n} \\
    \end{bmatrix}
    \in \R^{d_n \times d_n},
\end{equation}
where all diagonal blocks $G_{n,n}$ are square.  In other words, the
matrix $G_i$ arises from $G_{i-1}$ by appending a block column
(and adjusting the size).  We aim at computing the sequence of
matrix exponentials
\begin{equation}
    \label{eq:exp_sequence}
    \exp(G_0),\; \exp(G_1),\; \exp(G_2),\; \dotsc .
\end{equation}

One could, of course, simply compute each of the exponentials~\eqref{eq:exp_sequence} individually using standard techniques (see~\cite{Moler2003} for an overview).
However, the sequence of matrix exponentials~\eqref{eq:exp_sequence}
inherits the nested structure from the matrices $G_n$ in~\eqref{eq:G_n},
i.e., $\exp(G_{n-1})$ is a leading principle submatrix of $\exp(G_n)$.
In effect only the last block column of $\exp(G_n)$ needs to be computed and the goal of this paper is to explain how this can be achieved in a numerically safe manner.

In the special case where the spectra of the diagonal blocks $G_{n,n}$
are separated, Parlett's method~\cite{Parlett1976} yields -- in principle -- an efficient computational scheme:  Compute $F_{0,0} \defby
\exp(G_{0,0})$ and $F_{1,1} \defby \exp(G_{1,1})$ separately, then the
missing (1,2) block of $\exp(G_1)$ is given as the unique solution $X$ to
the Sylvester equation
\begin{equation*}
    G_{0,0} X - X G_{1,1} = F_{0,0} G_{0,1} - G_{0,1} F_{1,1}.
\end{equation*}
Continuing in this manner all the off-diagonal blocks required to compute~\eqref{eq:exp_sequence} could be
obtained from solving Sylvester equations. However, it is well known (see chapter 9 in~\cite{Higham2008}) that Parlett's method is numerically safe only when the spectra of the diagonal blocks are well separated, in the sense that all involved Sylvester equations are well conditioned. Since we consider the block structure as fixed, imposing such a condition would severely limit the scope of applications; it is certainly not met by the application we discuss below. 

A general class of applications for the described incremental
computation of exponentials arises from the matrix representations of
a linear operator $\mathcal{G} : V \to V$ restricted to a sequence of nested, finite
dimensional subspaces of a given infinite dimensional vector space $V$.
More precisely, one starts with a
finite dimensional subspace $V_0$ of $V$ with a basis $\mathcal{B}_0$.
Successively, the vector space $V_0$ is extended to $V_1 \subseteq V_2
\subseteq \cdots \subseteq V$ by generating a sequence of nested bases
$\mathcal{B}_0 \subseteq \mathcal{B}_1 \subseteq
\mathcal{B}_2\subseteq \cdots$. Assume that $\mathcal{G} V_n \subseteq V_n$ for all
$n=0,1,\dotsc$, and consider the sequence of matrix representations
$G_{n}$ of $\mathcal{G}$ with respect to $\mathcal{B}_n$.  Due to the
nestedness of the bases, $G_n$ is constructed from
$G_{n-1}$ by adding the columns representing the action of
$\mathcal{G}$ to $\mathcal{B}_n \setminus \mathcal{B}_{n-1}$.
As a
result, we obtain a sequence of matrices structured as
in~\eqref{eq:G_n}. 

A specific example for the scenario outlined above arises in computational finance, when pricing options
based on polynomial diffusion models; see~\cite{filipovic2016polynomial}. As we explain in more detail in section~\ref{sct3}, in this setting $\mathcal{G}$ is the generator of a stochastic differential equation (SDE), and $V_n$ are nested subspaces of multivariate polynomials.
Some pricing techniques require the computation of certain conditional moments that can be extracted from the matrix exponentials~\eqref{eq:exp_sequence}.
While increasing $n$ allows for a better approximation of the option price, the value of $n$ required to attain a desired accuracy is usually not known a priori. Algorithms that choose $n$ adaptively can be expected to rely on the incremental computation  of the whole sequence~\eqref{eq:exp_sequence}. 

Exponentials of block triangular matrices have also been studied in other contexts.  For two-by-two block triangular matrices, Dieci and
Papini study conditioning issues in~\cite{Dieci2001}, and a discussion
on the choice of scaling parameters for using Pad\'e approximants to
exponential function in~\cite{Dieci2000}.  In the case where the
matrix is also block-Toeplitz, a fast exponentiation algorithm is
developed in~\cite{Bini2016}.

The rest of this paper is organized as follows. In
section~\ref{sec:scaling_and_squaring} we give a detailed
description of our algorithm for incrementally computing exponentials
of block triangular matrices as in~\eqref{eq:G_n}. In section 3 we
discuss polynomial diffusion models, and some pricing
techniques which necessitate the use of such an incremental algorithm. Finally, 
numerical results are presented in section 4.

%% file: chapter2.tex
\section{Incremental scaling and squaring}
\label{sec:scaling_and_squaring}

Since the set of conformally partitioned block triangular matrices
forms an algebra, and $\exp(G_n)$ is a polynomial in $G_n$, the matrix
$\exp(G_n)$ has the same block upper triangular structure as $G_n$, that is,
\begin{equation*}
    \exp(G_n) =
    \begin{bmatrix}
        \exp(G_{0,0}) & \ast          & \cdots & \ast \\
                      & \exp(G_{1,1}) & \ddots & \vdots \\
                      &               & \ddots & \ast  \\
                      &               &        & \exp(G_{n,n}) \\
    \end{bmatrix}
    \in \R^{d_n \times d_n}.
\end{equation*}
As outlined in the introduction, we aim at computing $\exp(G_n)$ block
column by block column, from left to right.  Our algorithm is based on
the scaling and squaring methodology, which we briefly summarize
next.

\subsection{Summary of the scaling and squaring method}
\label{sec:scaling_squaring}

The scaling and squaring method uses a rational function to
approximate the exponential function, and typically involves three
steps. Denote by $r_{k,m}(z) = \frac{p_{k,m}(z)}{q_{k,m}(z)}$ the
$(k,m)$-Pad\'e approximant to the exponential function, meaning that
the numerator is a polynomial of degree $k$, and the denominator is a
polynomial of degree $m$. These Pad\'e approximants are very accurate close
to the origin, and in a first step the input matrix $G$ is therefore scaled by a
power of two, so that $\norm{2^{-s}G}$ is small enough to guarantee an
accurate approximation $r_{k,m}(2^{-s}G) \approx \exp(2^{-s} G)$.

The second step consists of evaluating the rational approximation
$r_{k,m}(2^{-s}G)$, and, finally, an approximation
to $\exp(G)$ is obtained in a third step by repeatedly squaring the
result, i.e.,
\begin{equation*}
    \exp(G) \approx r_{k,m}(2^{-s}G)^{2^s}.
\end{equation*}

Different choices of the scaling parameter $s$, and of the
approximation degrees $k$ and $m$ yield methods of different
characteristics. The choice of these parameters is critical for
the approximation quality, and for the computational efficiency,
see~\cite[chapter~10]{Higham2008}.

In what follows we describe techniques that allow for an incremental
evaluation of the matrix exponential of the block triangular
matrix~\eqref{eq:G_n}, using scaling and squaring.  These techniques
can be used with any choice for the actual underlying scaling and
squaring method, defined through the parameters $s$, $k$, and $m$.


\subsection{Tools for the incremental computation of exponentials}
\label{sec:tools}
Before explaining the algorithm, we first introduce some notation that
is used throughout.  The matrix $G_n$ from~\eqref{eq:G_n} can be written as
\begin{equation}
    \label{eq:G_twobytwo}
    G_n =
\left[
    \begin{array}{ccc|c}
        G_{0,0} &  \cdots &G_{0,n-1}& G_{0,n} \\
                &  \ddots & \vdots & \vdots  \\ 
                &   & G_{n-1,n-1}& G_{n-1,n} \\ \hline
                &         & & G_{n,n} \\
    \end{array}
\right]
    \bydef
    \begin{bmatrix}
        G_{n-1} & g_n \\
        0       & G_{n,n}
    \end{bmatrix}
\end{equation}
where $G_{n-1} \in \mathbb{R}^{d_{n-1} \times d_{n-1}}$, $G_{n,n}\in
\mathbb{R}^{b_n \times b_n}$, so that $g_n \in \R^{d_{n-1} \times
b_n}$.  Let $s$ be the scaling parameter, and $r = \frac{p}{q}$ the
rational function used in the approximation (for simplicity we will often omit the indices $k$ and $m$).  We denote the scaled
matrix by $\tilde{G}_n \defby 2^{-s} G_n$ and we partition it as
in~\eqref{eq:G_twobytwo}.\

The starting point of the algorithm consists in computing the Pad\'e approximant of the exponential $\exp(G_0)= \exp(G_{0,0})$, using a scaling and squaring method. Then, the sequence of matrix exponentials~\eqref{eq:exp_sequence} is incrementally computed by reusing at each step previously obtained quantities. So more generally, assume that $\exp(G_{n-1})$ has been approximated by using a scaling and squaring
method.  The three main computational steps for obtaining the Pad\'e approximant of $\exp(G_n)$ are
\begin{inparaenum}[(i)]
\item evaluating the polynomials $p(\tilde{G}_n)$, $q(\tilde{G}_n)$,
\item evaluating $p(\tilde{G}_n)^{-1} q(\tilde{G}_n)$, and
\item repeatedly squaring it.
\end{inparaenum}
We now discuss each of these steps separately, noting the quantities to keep at every iteration.

\subsubsection{Evaluating $p(\tilde{G}_n)$, $q(\tilde{G}_n)$ from $p(\tilde{G}_{n-1})$, $q(\tilde{G}_{n-1})$}
Similarly to \eqref{eq:G_twobytwo}, we start by writing $P_n\defby p(\tilde{G}_n)$ and $Q_n\defby q(\tilde{G}_n)$ as
\begin{equation*}
    P_n = 
    \begin{bmatrix}
        P_{n-1} & p_n \\
        0       & P_{n,n}
    \end{bmatrix},
    \quad
    Q_n = 
    \begin{bmatrix}
        Q_{n-1} & q_n \\
        0       & Q_{n,n}
    \end{bmatrix}.
\end{equation*}
In order to evaluate $P_{n}$, we first need to compute monomials of $\tilde{G}_{n}$, which for $l = 1, \dotsc, k$, can be written as
\begin{equation*}
    \tilde{G}_{n}^l =
    \begin{bmatrix}
        \tilde{G}_{n-1}^l & \sum_{j=0}^{l-1}
            \tilde{G}_{n-1}^j \tilde{g}_{n} \tilde{G}_{n,n}^{l-j-1} \\
        &                   \tilde{G}_{n,n}^l
    \end{bmatrix}.
\end{equation*}
Denote by $X_l \defby \sum_{j=0}^{l-1} \tilde{G}_{n-1}^j
\tilde{g}_{n} \tilde{G}_{n,n}^{l-j-1}$ the upper off diagonal block of
$\tilde{G}_n^l$, then we have the relation
\begin{equation*}
    X_{l} = \tilde{G}_{n-1}X_{l-1} + \tilde{g}_n \tilde{G}_{n,n}^{l-1}, \quad \text{for } l = 2, \cdots, k,
\end{equation*}
with $X_{1} \defby \tilde{g}_n$, so that all the monomials $\tilde{G}_n^l$, $l = 1, \dotsc, k$, can be
computed in $\bigO(b_n^3 + d_{n-1} b_n^2 + d_{n-1}^2 b_n)$.  Let $p(z)
= \sum_{l=0}^k \alpha_l z^l$ be the numerator polynomial of $r$, then
we have that
\begin{equation}
    \label{eq:P_n}
    P_n = 
    \begin{bmatrix}
        P_{n-1} & \sum_{l=0}^{k} \alpha_l X_l\\
                & p(\tilde{G}_{n,n})
    \end{bmatrix},
\end{equation}
which can be assembled in $\bigO(b_n^2 + d_{n-1} b_n)$, since only the
last block column needs to be computed. The complete evaluation of $P_n$ is summarized in Algorithm~\ref{alg:P_n}.
\begin{algorithm}[ht]
    \caption{Evaluation of $P_n$, using $P_{n-1}$
    \label{alg:P_n}}
    \begin{algorithmic}[1]
        \REQUIRE $G_{n-1}, G_{n,n}, g_n, P_{n-1}$, Pad\'e coefficients $\alpha_l, l = 0, \cdots, k$.
        \ENSURE $P_n$.
        \STATE $\tilde{g}_n \leftarrow 2^{-s} g_n$,
            $\tilde{G}_{n,n} \leftarrow 2^{-s} G_{n,n}$,
	 $\tilde{G}_{n-1} \leftarrow 2^{-s} G_{n-1}$
        \STATE $X_1 \leftarrow \tilde{g}_n$
        \FOR {$l=2, 3,  \cdots, k$}
	\STATE Compute $\tilde{G}_{n,n}^l$
         	\STATE $X_l = \tilde{G}_{n-1} X_{l-1} + \tilde{g}_n \tilde{G}_{n,n}^{l-1}$
        \ENDFOR
        \STATE $X_0 \leftarrow \mathbf{0}_{d_{n-1} \times b_n}$
        \STATE Compute off diagonal block of $P_n$:  $\sum_{l=0}^{k} \alpha_l X_l$.
        \STATE Compute $p(\tilde{G}_{n,n}) = \sum_{l=0}^{k} \alpha_l \tilde{G}_{n,n}^l$
        \STATE Assemble $P_n$ as in $\eqref{eq:P_n}$
    \end{algorithmic}
\end{algorithm}

Similarly, one computes $Q_n$ from $Q_{n-1}$, using again the matrices
$X_l$.
\subsubsection{Evaluating $Q_n^{-1} P_n$}

With the matrices $P_n$, $Q_n$ at hand, we now need to compute the
rational approximation $Q_n^{-1} P_n$. 
We assume that $Q_n$ is well
conditioned, in particular non-singular, which is ensured by the choice of the 
scaling parameter and of the Pad\'e approximation, see, e.g.,~\cite{Higham2009}.
We focus on the computational cost.
For simplicity, we introduce the notation
\begin{equation*}
    \tilde{F}_n = 
    \begin{bmatrix}
        \tilde{F}_{0,0} &  \cdots & \tilde{F}_{0,n} \\
                        &  \ddots & \vdots  \\
                        &         & \tilde{F}_{n,n} \\
    \end{bmatrix}
    \defby Q_n^{-1} P_n,
    \quad
    F_n = 
    \begin{bmatrix}
        F_{0,0} &  \cdots & F_{0,n} \\
                &  \ddots & \vdots  \\
                &         & F_{n,n} \\
    \end{bmatrix}
    \defby \tilde{F}_n^{2^s},
\end{equation*}
and we see that
\begin{equation}
    \label{eq:F_n_tilde}
    \begin{split}
    \tilde{F}_n & = Q_n^{-1} P_n =
    \begin{bmatrix}
        Q_{n-1}^{-1} & - Q_{n-1}^{-1} q_n Q_{n,n}^{-1} \\
        0            & Q_{n,n}^{-1}
    \end{bmatrix}
    \begin{bmatrix}
        P_{n-1} & p_n \\
        0       & P_{n,n}
    \end{bmatrix}\\
    & =
    \begin{bmatrix}
        \tilde{F}_{n-1} & Q_{n-1}^{-1} ( p_n - q_n Q_{n,n}^{-1} P_{n,n} ) \\
        0               & Q_{n,n}^{-1} P_{n,n}
    \end{bmatrix}.
    \end{split}
\end{equation}

To solve the linear system $Q_{n,n}^{-1} P_{n,n}$ we compute an LU
decomposition with partial pivoting for
$Q_{n,n}$, requiring $\bigO(b_n^3)$ operations.
This LU decomposition is saved for future use, and hence we may assume
that we have available the LU decompositions for all diagonal
blocks from previous computations:
\begin{equation}
    \label{eq:store_lu}
    \Pi_l Q_{l,l}= L_l U_l, \quad l=0, \dotsc, n-1.
\end{equation}
Here, $\Pi_l \in \mathbb{R}^{b_l \times b_l},$ $l=0, \dotsc, n-1$ are permutation matrices; 
$L_l\in \mathbb{R}^{b_l \times b_l},$ $l=0, \dotsc, n-1$ are lower triangular matrices and $U_l\in \mathbb{R}^{b_l \times b_l},$ $l=0, \dotsc, n-1$ are upper triangular matrices.

Set $Y_n \defby p_n - q_n Q_{n,n}^{-1} P_{n,n} \in \R^{d_{n-1} \times
b_n}$, and partition it as
\begin{equation*}
    Y_n = 
    \begin{bmatrix}
      Y_{0,n} \\
      \vdots  \\
      Y_{n-1,n} \\
    \end{bmatrix}.
\end{equation*}
Then we compute $Q_{n-1}^{-1} Y_n$ by block backward
substitution, using the decompositions of the diagonal blocks.  The
total number of operations for this computation is hence $\bigO(d_{n-1}^2 b_n + d_{n-1} b_n^2)$, so that the number of
operations for computing $\tilde{F}_n$ is $\bigO(b_n^3 + d_{n-1}^2
b_n + d_{n-1} b_n^2)$.
Algorithm~\ref{alg:tilde_F_n} describes the complete procedure to compute $\tilde{F}_n$.

\begin{algorithm}[ht]
    \caption{Evaluation of $\tilde{F}_n = Q_n^{-1} P_n$}
    \label{alg:tilde_F_n}
    \begin{algorithmic}[1]
        \REQUIRE $Q_n, P_n$ and quantities \eqref{eq:store_lu}
        \ENSURE $\tilde{F}_n = Q_n^{-1} P_n$ and LU decomposition of $Q_{n,n}$.
        \STATE Compute $\Pi_n Q_{n,n}= L_n U_n$ and keep it for future use \eqref{eq:store_lu}
        \STATE Compute  $\tilde{F}_{n,n} \defby Q_{n,n}^{-1} P_{n,n}$
        \STATE $Y_n = p_n - q_n Q_{n,n}^{-1} P_{n,n}$
        \STATE $\tilde{F}_{n-1,n} = U_{n-1}^{-1} L_{n-1}^{-1} \Pi_{n-1} Y_{n-1,n} $
        \FOR {$l = n - 2, n - 3,  \cdots, 0$}
        \STATE $\tilde{F}_{l, n} = U_{l}^{-1}L_{l}^{-1} \Pi_{l}(Y_{l, n}-\sum_{j=l+1}^{n-1} Q_{l,j}\tilde{F}_{j, n}) $
        \ENDFOR
        \STATE Assemble $\tilde{F}_n$ as in \eqref{eq:F_n_tilde}
    \end{algorithmic}
\end{algorithm}
\subsubsection{The squaring phase}

Having computed $\tilde{F}_n$, which we write as
\begin{equation*}
    \tilde{F}_n =
    \begin{bmatrix}
        \tilde{F}_{n-1} & \tilde{f}_n \\
                        & \tilde{F}_{n,n}
    \end{bmatrix},
\end{equation*}
we now need to compute $s$ repeated squares of that matrix, i.e.,
\begin{equation}
    \label{eq:squares}
    \tilde{F}_{n}^{2^l} =
    \begin{bmatrix}
        \tilde{F}_{n-1}^{2^l} & \sum_{j=0}^{l-1}
            \tilde{F}_{n-1}^{2^{l - 1 + j}} \tilde{f}_n \tilde{F}_{n,n}^{2^j}\\
                              & \tilde{F}_{n,n}^{2^l}
   \end{bmatrix}, \quad l = 1, \dotsc, s,
\end{equation}
so that $F_n = \tilde{F}_{n}^{2^s}$.  Setting $Z_l \defby
\sum_{j=0}^{l-1} \tilde{F}_{n-1}^{2^{l - 1 + j}} \tilde{f}_j
\tilde{F}_{n,n}^{2^j}$, we have the recurrence
\begin{equation*}
    Z_l = \tilde{F}_{n-1}^{2^{l-1}} Z_{l-1} + Z_{l-1} \tilde{F}_{n,n}^{2^{l-1}},
\end{equation*}
with $Z_0 \defby \tilde{f}_n$.  Hence, if we have stored the
intermediate squares from the computation of $F_{n-1}$, i.e.,
\begin{equation}
    \label{eq:store_squares}
    \tilde{F}_{n-1}^{2^l},  \quad l=1, \dotsc, s
\end{equation}
we can compute all the quantities $Z_l$, $l=1, \dotsc, s$ in
$\bigO(d_{n-1}^2 b_n + d_{n-1} b_n^2)$, so that the total cost for
computing $F_n$ (and the intermediate squares of $\tilde{F}_n$) is
$\bigO(d_{n-1}^2 b_n + d_{n-1} b_n^2 + b_n^3)$. Again, we summarize the squaring phase in the following algorithm.

\begin{algorithm}[ht]
    \caption{Evaluation of $F_n = \tilde{F}_n^{2^s}$
    \label{alg:F_n}}
    \begin{algorithmic}[1]
        \REQUIRE $\tilde{F}_{n-1}, \tilde{f}_n, \tilde{F}_{n,n}$, quantities \eqref{eq:store_squares}.
        \ENSURE $F_n$ and updated intermediates.
        \STATE $Z_0 \leftarrow \tilde{f}_n$
        \FOR {$l=1, 2,  \cdots, s$}
          \STATE Compute $\tilde{F}_{n,n}^{2^l}$
          \STATE $Z_l = \tilde{F}_{n-1}^{2^{l-1}} Z_{l-1} + Z_{l-1} \tilde{F}_{n,n}^{2^{l-1}}$
          \STATE Assemble $\tilde{F}_n^{2^l}$ as in \eqref{eq:squares} and save it
        \ENDFOR
          \STATE $F_n \leftarrow \tilde{F}_{n}^{2^s}$
    \end{algorithmic}
\end{algorithm}
\subsection{Overall Algorithm}
\label{sec:overall}

Using the techniques from the previous section, we now give a concise
description of the overall algorithm.  We assume that the quantities
listed in equations~\eqref{eq:store_lu}
and~\eqref{eq:store_squares} are stored in memory, with a space requirement of $\bigO(d_{n-1}^2)$. 

In view of this, we assume that $F_{n-1}$ and the aforementioned intermediate
quantities have been computed.  Algorithm~\ref{alg:step} describes the
overall procedure to compute $F_n$, and to update the intermediates;
we continue to use the notation introduced in~\eqref{eq:G_twobytwo}.

\begin{algorithm}[ht]
    \caption{Computation of $F_n \approx \exp(G_n)$, using $F_{n-1}$
    \label{alg:step}}
    \begin{algorithmic}[1]
        \REQUIRE Block column $g_n$, diagonal block $G_{n,n}$, quantities \eqref{eq:store_lu},
        and~\eqref{eq:store_squares}.
        \ENSURE $F_n$, and updated intermediates.
        \STATE Extend $P_{n-1}$ to $P_n$ using Algorithm~\ref{alg:P_n},
        and form analogously $Q_n$
        \STATE Compute $\tilde{F}_n$ using Algorithm~\ref{alg:tilde_F_n}
        \STATE Evaluate $F_n = \tilde{F}_n^{2^s}$
        using Algorithm~\ref{alg:F_n}
    \end{algorithmic}
\end{algorithm}

As explained in the previous section, the number of operations for
each step in Algorithm~\ref{alg:step} is $\bigO(d_{n-1}^2 b_n +
d_{n-1} b_n^2 + b_n^3)$, using the notation at the beginning of
section~\ref{sec:tools}.  If $F_n$ were simply computed from scratch, without
the use of the intermediates, the number of operations for scaling and
squaring would be $\bigO((d_{n-1} + b_n)^3)$.  In the typical
situation where $d_{n-1} \gg b_n$, the dominant term in the latter
complexity bound is $d_{n-1}^3$, which is absent from the complexity
bound of Algorithm~\ref{alg:step}.

In order to solve our original problem, the computation of the
sequence $\exp(G_0)$, $\exp(G_1)$, $\exp(G_2)$, $\dotsc$, we use
Algorithm~\ref{alg:step} repeatedly; the resulting procedure is shown
in Algorithm~\ref{alg:full}.

\begin{algorithm}[ht]
    \caption{Approximation of $\exp(G_0), \exp(G_1), \dotsc$
    \label{alg:full}}
    \begin{algorithmic}[1]
        \REQUIRE Pad\'e approximation parameters $k$, $m$, and $s$
        \ENSURE $F_0 \approx \exp(G_0)$, $F_1 \approx \exp(G_1), \dotsc$
        \STATE Compute $F_0$ using scaling and squaring, store
            intermediates for Algorithm~\ref{alg:step}
        \FOR{$n=1,2,\dotsc$}
            \STATE Compute $F_{n}$ from $F_{n-1}$ using Algorithm~\ref{alg:step}
            \IF{termination criterion is satisfied} \RETURN
            \ENDIF
        \ENDFOR
    \end{algorithmic}
\end{algorithm}

We now derive a complexity bound for the number of operations
spent in Algorithm~\ref{alg:full}.  For simplicity of notation we
consider the case where all diagonal blocks are of equal size, i.e.,
$b_k \equiv b \in \N$, so that $d_k = (k+1)b$.  At iteration $k$ the
number of operations spent within Algorithm~\ref{alg:step} is thus 
$\bigO(k^2 b^3)$.  Assume that the termination criterion used in
Algorithm~\ref{alg:full} effects to stop the procedure after the
computation of $F_n$.  The overall complexity bound for the number of
operations until termination is $\bigO(\sum_{k=0}^n k^2 b^3 ) =
\bigO(n^3 b^3)$, which matches the complexity bound of applying
scaling and squaring only to $G_n \in \R^{(n+1)b \times (n+1)b}$,
which is also $\bigO( (nb)^3 )$.

In summary the number of operations needed to compute $F_n$ by
Algorithm~\ref{alg:full} is asymptotically the same as applying the
same scaling and squaring setting \emph{only} to compute $\exp(G_n)$,
while Algorithm~\ref{alg:full} incrementally reveals \emph{all}
exponentials $\exp(G_0)$, $\dotsc$, $\exp(G_n)$ in the course of the
iteration, satisfying our requirements outlined in the introduction.

\subsection{Adaptive scaling}

In Algorithms~\ref{alg:step} and~\ref{alg:full} we have assumed
that the scaling power $s$ is given as input parameter, and that it is
fixed throughout the computation of $\exp(G_0), \dotsc, \exp(G_n)$.
This is in contrast to what is usually intented in the scaling and squaring method, see Section~\ref{sec:scaling_squaring}.  On the one hand $s$
must be sufficiently large so that $r_{k,m}(2^{-s} G_l) \approx \exp(2^{-s}
G_l)$, for $0 \le l \le n$.  If, on the other hand, $s$ is chosen
\emph{too large}, then the evaluation of $r_{k,m}(2^{-s}G_l)$ may become
inaccurate, due to \emph{overscaling}.
So if $s$ is fixed, and the norms $\norm{G_l}$ grow with
increasing $l$, as one would normally expect, an accurate approximation cannot
be guaranteed for all $l$.

Most scaling and squaring designs hence choose $s$ in dependence of the
norm of the input matrix~\cite{Moler2003,Guttel2016,Higham2009}.  For
example, in the algorithm of Higham described in~\cite{Higham2009}, it
is the smallest integer satisfying
\begin{equation}
    \label{eqn:theta}
    \norm{2^{-s} G_l}_1 \le \theta \approx 5.37... .
\end{equation}
In order to combine our incremental evaluation techniques with this
scaling and squaring design, the scaling power $s$ must thus be chosen
dynamically in the course of the evaluation.  Assume that $s$
satisfies the criterion~\eqref{eqn:theta} at step $l-1$, but not at step
$l$.  We then simply discard all accumulated data structures from
Algorithm~\ref{alg:step}, increase $s$ to match the
bound~\eqref{eqn:theta} for $G_{l}$, and start
Algorithm~\ref{alg:full} anew with the \emph{repartitioned} input matrix
\begin{equation}
    \label{eqn:repart}
G_n = \left[
    \begin{array}{ccc|c|c|c}
        G_{0,0} &  \cdots &G_{0,l}  & G_{0,l+1}  & \cdots & G_{0,n}\\
                &  \ddots & \vdots  & \vdots     &        & \vdots\\ 
                &         & G_{l,l} & G_{l,l+1}  & \cdots & G_{l,n}\\ \hline
                &         &         & G_{l+1,l+1}& \cdots & G_{l+1,n}\\
                &         &         &            & \ddots & \vdots   \\
                &         &         &            &        & G_{n,n}  \\
    \end{array}
\right]
    =
    \underbrace{
\left[
    \begin{array}{c|c|c|c}
        \hat{G}_{0,0} & \hat{G}_{0,1}    & \cdots & \hat{G}_{0,n-l}\\ \hline
                      & \hat{G}_{1,1}    & \cdots & \hat{G}_{1,n-l}\\
                      &                  & \ddots & \vdots   \\
                      &                  &        & \hat{G}_{n-l,n-l}  \\
    \end{array}
\right]}_{\bydef \hat{G}_{n-l}}.
\end{equation}
The procedure is summarized in Algorithm~\ref{alg:adaptive}.

\begin{algorithm}[t]
    \caption{Approximation of $\exp(G_0), \exp(G_1), \dotsc$ with
    adaptive scaling
    \label{alg:adaptive}}
    \begin{algorithmic}[1]
        \REQUIRE Pad\'e approximation parameters $k$, $m$, norm
        bound $\theta$.
        \ENSURE $F_0 \approx \exp(G_0)$, $F_1 \approx \exp(G_1), \dotsc$
        \STATE $s \leftarrow \max\{0, \log(\norm{G_0}_1)\}$
        \STATE Compute $F_0$ using scaling and squaring, store
            intermediates for Algorithm~\ref{alg:step}
        \FOR{$l=1,2,\dotsc$}
            \IF{ $\norm{G_l}_1 > \theta$ }
                \STATE Repartition $G_n = \hat{G}_{n-l}$ as in~\eqref{eqn:repart}
                \STATE Restart algorithm with $\hat{G}_{n-l}$.
            \ENDIF
            \STATE Compute $F_{l}$ from $F_{l-1}$ using Algorithm~\ref{alg:step}
            \IF{termination criterion is satisfied} \RETURN
            \ENDIF
        \ENDFOR
    \end{algorithmic}
\end{algorithm}

It turns out that the computational overhead induced by this
restarting procedure is quite modest.  In the notation introduced for
the complexity discussion in Section~\ref{sec:overall}, the number of
operations for computing $\exp(G_n)$ by Higham's scaling and squaring
method is $\bigO(\log(\norm{G_n}_1) (nb)^3)$.  Since there are at
most $\log(\norm{G_n}_1)$ restarts in Algorithm~\ref{alg:adaptive},
the total number of operations for incrementally computing all
exponentials $\exp(G_0), \dotsc, \exp(G_n)$ can be bounded by a
function in $\bigO(\log(\norm{G_n}_1)^2 (nb)^3)$.  We assess the
actual performance of Algorithm~\ref{alg:adaptive} in
Section~\ref{sec:numerical_experiments}.

In our application from option pricing it turns out that the norms of
the matrices $G_l$ do not grow dramatically (see
Sections~\ref{sec:jacobix} and~\ref{sec:hestonx}) and quite accurate
approximations to all the matrix exponentials can be computed even if
the scaling factor is fixed (see Section~\ref{sec:exp_jacobi}).

%% file: chapter3.tex
\section{Option pricing in polynomial models}\label{sct3}

The main purpose of this section is to explain how certain option pricing techniques require the sequential computation of matrix exponentials for block triangular matrices. The description will necessarily be rather brief; we refer to, e.g., the textbook~\cite{Elliot2005} for more details.

Because we are evaluating at initial time $t=0$, the price of a certain option expiring at time $\tau>0$ consists of computing an expression of the form
\begin{equation}\label{priceexpectation}
e^{-r\tau} \mathbb{E}[f(X_\tau)],
\end{equation}
where $(X)_{0 \leq t \leq \tau}$ is a $d$-dimensional stochastic process modelling the price of financial assets over the time interval $[0,\tau]$, $f : \mathbb{R}^d \to \mathbb{R}$ is the so-called payoff function and $r$ represents a fixed interest rate.  In the following, we consider stochastic processes described by an SDE of the form
\begin{align}\label{SDE}
dX_t=b(X_t)dt+\Sigma(X_t)dW_t,
\end{align}
where $W$ denotes a $d$-dimensional Brownian motion, $b : \mathbb{R}^d \mapsto \mathbb{R}^{d}$, and $\Sigma : \mathbb{R}^d \mapsto \mathbb{R}^{d \times d}$.

\subsection{Polynomial diffusion models}

During the last years, polynomial diffusion models have become a versatile tool in financial applications, including option pricing. In the following, we provide a short summary and refer to the paper by Filipovi\'c and Larsson \cite{filipovic2016polynomial} for the mathematical foundations. 

For a polynomial diffusion process one assumes that the coefficients of the vector $b$ in~\eqref{SDE} and the matrix $A : = \Sigma \Sigma^T$ satisfy
\begin{equation}\label{polynomial}
A_{ij} \in \text{Pol}_2(\mathbb{R}^d), \qquad b_i \in \text{Pol}_1(\mathbb{R}^d)  \quad \text{for} \quad i,j = 1,\ldots,d.
\end{equation}
Here, $\text{Pol}_n(\mathbb{R}^d)$ represents the set of $d$-variate polynomials of total degree at most $n$, that is,
\begin{equation*}
\text{Pol}_n(\mathbb{R}^d):= \Big\{\sum_{0 \le |\textbf{k}|\le n} \alpha_{\textbf{k}} x^{\bf{k}}| x \in \mathbb{R}^d, \alpha_{\textbf{k}} \in \mathbb{R}\Big\},
\end{equation*}
where we use multi-index notation: $\mathbf{k}=(k_1, \dots, k_d) \in \mathbb{N}_0^d$, $|\mathbf{k}|:=k_1+\dots+k_d$ and $x^{\bf{k}}:=x_1^{k_1}\dots x_d^{k_d}$. In the following, $\text{Pol}(\mathbb{R}^d)$ represents the set of all multivariate polynomials on $\mathbb{R}^d$.\

Associated with $A$ and $b$ we define the partial differential operator $\mathcal{G}$ by 
\begin{equation}\label{generator}
\mathcal{G}f=\frac{1}{2}\tr(A \nabla^2 f)+b^T \nabla f.
\end{equation}
which represents the so called generator
for~\eqref{SDE}, see~\cite{Oksendal2003}. It can be directly verified that~\eqref{polynomial} implies that $\text{Pol}_n(\mathbb{R}^d$) is invariant under $\mathcal{G}$  for any $n \in \mathbb{N}$, that is, 
\begin{equation}\label{PreservingProperty}
\mathcal{G}\text{Pol}_n(\mathbb{R}^d) \subseteq \text{Pol}_n(\mathbb{R}^d).
\end{equation}
\begin{remark}
    \begin{rm}
In many applications, one is interested in solutions to $\eqref{SDE}$ that lie on a state space $E \subseteq \mathbb{R}^d$ to incorporate, for example, nonnegativity. This problem is largely studied in~\cite{filipovic2016polynomial}, where existence and uniqueness of solutions to \eqref{SDE} on several types of state spaces $E \subseteq \mathbb{R}^d$ and for large classes of $A$ and $b$ is shown.
    \end{rm}
\end{remark}
Let us now fix a basis of polynomials $\mathcal{H}_n=\{h_1, \dots, h_N\}$ for $\text{Pol}_n(\mathbb{R}^d)$, where $N= \dim \text{Pol}_n(\mathbb{R}^d) = \binom{n+d}{n}$, and write
\begin{equation*}
H_n(x)=(h_1(x), \dots, h_N(x))^T.
\end{equation*}
Let $G_n$ denote the matrix representation with respect to $\mathcal{H}$ of the linear operator $\mathcal{G}$ restricted to $\text{Pol}_n(\mathbb{R}^d)$. By definition,
\begin{equation*}
\mathcal{G}p(x)=H_n(x)^T G_n \vec{p}.
\end{equation*}
for any $p\in \text{Pol}_n(\mathbb{R}^d)$ with coordinate vector $\vec{p} \in \R^N$ with respect to $\mathcal{H}_n$.
By Theorem 3.1 in \cite{filipovic2016polynomial}, the  corresponding polynomial moment can be computed from
\begin{equation}\label{condmoments}
\mathbb{E}[p(X_\tau)]=H_n(X_0)^Te^{\tau G_n} \vec{p}.
\end{equation}

The setting discussed above corresponds to the scenario described in the introduction. We have a sequence of subspaces
\[ \text{Pol}_0(\mathbb{R}^d) \subseteq \text{Pol}_1(\mathbb{R}^d) \subseteq \text{Pol}_2(\mathbb{R}^d) \subseteq \cdots \subseteq \text{Pol}(\mathbb{R}^d)\]
and the polynomial preserving property~\eqref{PreservingProperty} implies that the matrix representation $G_n$ is block upper triangular with $n+1$ square diagonal blocks of size \[1, d, \binom{1+d}{2}, \ldots, \binom{n+d-1}{n}.\]

In the rest of this section we introduce two different pricing techniques that require the incremental computation of polynomial moments of the form~\eqref{condmoments}. 

\subsection{Moment-based option pricing for Jacobi models} \label{sec:jacobix}
The Jacobi stochastic volatility model is a special case of a polynomial diffusion model and it is characterized by the SDE
\begin{align*}
&dY_t=(r-V_t/2)dt + \rho \sqrt{Q(V_t)} dW_{1t}+\sqrt{V_t-\rho^2Q(V_t)}dW_{2t},\\
&dV_t=\kappa(\theta -V_t)dt+\sigma \sqrt{Q(V_t)}dW_{1t},
\end{align*}
where 
\begin{equation*}
Q(v)=\frac{(v-v_{\min})(v_{\max}-v)}{(\sqrt{v_{\max}}-\sqrt{v_{\min}})^2},
\end{equation*}
for some $0 \leq v_{\min} < v_{\max}$. Here, $W_{1t}$ and $W_{2t}$ are independent standard Brownian motions and the model parameters satisfy the conditions $\kappa \geq 0$, $\theta \in [v_{\min},v_{\max}]$, $\sigma >0$, $r \geq 0$, $\rho \in [-1,1]$.\
In their paper, Ackerer et al.~\cite{ackerer2016jacobi} use this model in the context of option pricing where the price of the asset is specified by $S_t \defby e^{Y_t}$ and $V_t$ represents the squared stochastic volatility. In the following, we briefly introduce the pricing technique they propose and  explain how it involves the incremental computation of polynomial moments. 

Under the Jacobi model with the discounted payoff function $f$ of an European claim, 
the option price~\eqref{priceexpectation} at initial time $t=0$ can be expressed as 
\begin{equation}\label{price}
\sum_{n \geq 0} f_n l_n,
\end{equation}
where $\{f_n, n\geq 0\}$ are the Fourier coefficients of $f$ and $\{l_n, n\geq 0\}$ are Hermite moments. As explained in~\cite{ackerer2016jacobi}, the Fourier coefficients can be conveniently computed in a recursive manner. The Hermite moments are computed using~\eqref{condmoments}. Specifically,
consider the monomial basis of $\text{Pol}_n(\R^2)$:
\begin{equation}\label{eq:basisJacobi}
H_n(y,v) \defby (1,y,v,y^2,yv,v^2,\dots,y^n,y^{n-1}v,\dots,v^n)^T.
\end{equation}
Then 
\begin{equation} \label{eq:hermitemoments}
l_n = H_n(Y_0,V_0)^Te^{\tau G_n} \vec{h}_n,
\end{equation}
where $\vec{h}_n$ contains the coordinates with respect to~\eqref{eq:basisJacobi} of 
\begin{equation*}
 \frac{1}{\sqrt{n!}} h_n \left(\frac{y-\mu_w}{\sigma_w}\right),
\end{equation*}
with real parameters $\sigma_w, \mu_w$ and the $n$th Hermite polynomial ${h}_n$.

Truncating the sum \eqref{price} after a finite number of terms allows
one to obtain an approximation of the option price.
Algorithm~\ref{alg:jacobi} describes a heuristic to selecting the truncation based on the absolute value of the summands, using Algorithm~\ref{alg:full} for computing the required moments incrementally. 
\begin{algorithm}[ht]
\caption{Option pricing for the European call option under the Jacobi stochastic volatility model}
\begin{algorithmic}[1]\label{alg:jacobi}
\REQUIRE Model and payoff parameters, tolerance $\epsilon$
\ENSURE Approximate option price
\STATE $n=0$
\STATE Compute $l_0$, $f_0$; set $\text{Price} = l_0 f_0$.
\WHILE{$|l_n f_n| > \epsilon \cdot \text{Price}$}
\STATE $n=n+1$
\STATE \label{line:jacobiexp} Compute $\exp(\tau G_n)$ using Algorithm~\ref{alg:step}.
\STATE \text{Compute Hermite moment $l_n$ using~\eqref{eq:hermitemoments}. }
\STATE \text{Compute Fourier coefficient $f_n$ as described in~\cite{ackerer2016jacobi}.}
\STATE $\text{Price}=\text{Price}+l_n f_n$;
\ENDWHILE
\end{algorithmic}
\end{algorithm}

As discussed in Section~\ref{sec:scaling_and_squaring}, a norm estimate for $G_n$ is instrumental for choosing a priori the scaling parameter in the scaling and squaring method. The following lemma provides such an estimate for the model under consideration. 
\begin{lemma} \label{lemmanormJ}
Let $G_n$ be the matrix representation of the operator $\mathcal{G}$ defined in~\eqref{generator}, with respect to the basis~\eqref{eq:basisJacobi} of $\mathrm{Pol}_{n}(\mathbb{R}^2)$.  Define 
\begin{equation*}
\alpha:=\frac{\sigma (1+v_{\min}v_{\max}+v_{\max}+v_{\min})}{2(\sqrt{v_{\max}}-\sqrt{v_{\min}})^2}.
\end{equation*}
Then the matrix 1-norm of $G_n$ is bounded by
\begin{equation*}
n( r + \kappa + \kappa\theta - \sigma \alpha ) + \frac12 n^2 ( 1 + |\rho| \alpha + 2 \sigma \alpha).
\end{equation*}
\begin{proof}
The operator $\mathcal{G}$ in the Jacobi model takes the form
\begin{equation*}
\mathcal{G}f(y,v)=\frac{1}{2}\tr(A(v) \nabla^2 f(y,v))+b(v)^\top \nabla f(y,v),
\end{equation*}
where 
\begin{equation*}
b(v)=\begin{bmatrix}
  r-v/2 \\
  \kappa(\theta -v)
\end{bmatrix},
 \quad 
A(v)=\begin{bmatrix}
  v & \rho \sigma Q(v)\\
 \rho \sigma Q(v) &  \sigma^2 Q(v)
\end{bmatrix}.
\end{equation*} 
Setting $S\defby (\sqrt{v_{\max}}-\sqrt{v_{\min}})^2$, we consider the action of the generator $\mathcal{G}$ on a basis element $y^pv^q$:
\begin{align*}
\mathcal{G} y^p v^q=&y^{p-2} v^{q+1}p \frac{p-1}{2} - y^{p-1}v^{q+1} p \Big(\frac{1}{2}+\frac{q\rho \sigma}{S}\Big)+y^{p-1} v^qp \Big(r+q \rho \sigma \frac{v_{\max} + v_{\min}}{S}\Big)\\
& -y^{p-1}v^{q-1}\frac{pq \rho \sigma v_{\max} v_{\min}}{S} - y^p v^q q \Big(\kappa + \frac{q-1}{2} \frac{\sigma^2}{S} \Big) \\
& - y^p v^{q-2} q \frac{q-1}{2} \frac{\sigma^2 v_{\max}v_{\min}}{S} +y^p v^{q-1} q \Big(\kappa \theta + \frac{q-1}{2} \sigma^2 \frac{ v_{\max} + v_{\min}}{S} \Big). 
\end{align*}
For the matrix 1-norm of $G_n$, one needs to determine the values of $(p,q) \in \mathcal{M}:= \{(p,q) \in \mathbb{N}_0 \times \mathbb{N}_0 | p+q \leq n\}$ for which the $1$-norm of the coordinate vector of $\mathcal{G} y^p v^q$ becomes maximal. Taking into account the nonnegativity of the involved model parameters and replacing $\rho$ by $|\rho|$, we obtain an upper bound as follows:
\begin{align*}
& p \frac{p-1}{2} + p \Big(\frac{1}{2}+\frac{q|\rho|\sigma}{S}\Big)+ p \Big(r+q |\rho|\sigma \frac{v_{\max} + v_{\min}}{S}\Big)+ \frac{pq |\rho|\sigma v_{\max} v_{\min}}{S}\\
&+ q \Big(\kappa + \frac{q-1}{2} \frac{\sigma^2}{S} \Big) + q \frac{q-1}{2} \frac{\sigma^2 v_{\max}v_{\min}}{S} + q \Big(\kappa \theta + \frac{q-1}{2} \sigma^2 \frac{ v_{\max} + v_{\min}}{S} \Big) \\
=& pr + q\kappa(\theta+1) + \frac12 p^2 + 2pq|\rho| \alpha +  q (q-1) \sigma \alpha  \\
\le & n( r + \kappa + \kappa\theta ) + \frac12 n^2 + 2pq|\rho| \alpha + n(n-1)\sigma \alpha. 
\end{align*}
This completes the proof, noting that the maximum of $pq$ on $\mathcal{M}$ is bounded by $n^2 / 4$
over $\mathcal{M}$.
\end{proof}
\end{lemma}

The result of Lemma~\ref{lemmanormJ} predicts that the norm of $G_n$ grows, in general, quadratically. This prediction is confirmed numerically for parameter settings of practical relevance.
\subsection{Moment-based option pricing for Heston models} \label{sec:hestonx}

The Heston model is another special case of a polynomial diffusion model, characterized by the SDE 
\begin{align*}
&dY_t=(r-V_t/2)d_t + \rho \sqrt{V_t} dW_{1t}+\sqrt{V_t}\sqrt{1-\rho^2}dW_{2t},\\
&dV_t=\kappa(\theta -V_t)dt+\sigma \sqrt{V_t}dW_{1t},
\end{align*}
with model parameters satisfying the conditions $\kappa \geq 0$, $\theta \geq 0$, $\sigma >0$, $r \geq 0$, $\rho \in [-1,1]$. As before, the asset price is modeled via $S_t \defby e^{Y_t}$, while $V_t$ represents the squared stochastic volatility.

Lasserre et al.~\cite{Lasserre2006} developed a general option pricing technique based on moments and semidefinite programming (SDP). In the following we briefly explain the main steps and in which context an incremental computation of moments is needed. In doing so, we restrict ourselves to the specific case of the Heston model and European call options.

Consider the payoff function $f(y) \defby (e^{y}-e^K)^+$ for a certain log strike value $K$. Let $\nu(dy)$ be the $Y_\tau$-marginal distribution of the joint distribution of the random variable $(Y_\tau,V_\tau)$. Define the restricted measures $\nu_1$ and $\nu_2$ as $\nu_1=\nu |_{(- \infty , K]}$ and $\nu_2=\nu |_{[K, \infty)}$. By approximating the exponential in the payoff function with a Taylor series truncated after $n$ terms, the option price~\eqref{priceexpectation} can be written as a certain linear function $L$ in the moments of $\nu_1$ and $\nu_2$, i.e., 
\begin{equation*}
\mathbb{E}[f(Y_\tau)]= L(n,\nu_1^0, \cdots, \nu_1^n,\nu_2^0, \cdots, \nu_2^n),
\end{equation*}
where $\nu_i^m$ represents the $m$th moment of the $i$th measure.\

A lower / upper bound of the option price can then be computed by solving the optimization problems 
 \begin{align} \label{optiproblem}
 SDP_n \defby \left\{
                \begin{array}{ll}
                  \min/ \max \hspace{0.1 cm}&L(n,\nu_1^0, \cdots, \nu_1^n,\nu_2^0, \cdots, \nu_2^n)\\
\text{subject to   }\hspace{0.1 cm}  &\nu_1^j+\nu_2^j=\nu^j, \quad j=0, \cdots, n\\ 
                  &\nu_1 \text{ is a Borel measure on } (- \infty , K],\\
                  &\nu_2 \text{ is a Borel measure on } [K,\infty).\\
		\end{array}
              \right.
\end{align}
Two SDP arise when writing the last two conditions in \eqref{optiproblem} via moment and localizing matrices, corresponding to the so-called truncated Stieltjes moment problem.

Formula \eqref{condmoments} is used in this setting to compute the moments $\nu^j$. Increasing the relaxation order $n$ iteratively allows us to find sharper bounds (this is trivial because increasing $n$ adds more constraints). One stops as soon as the bounds are sufficiently close. Algorithm~\ref{algoSDP} summarizes the resulting pricing algorithm.
\begin{algorithm}[ht]
\caption{Option pricing for European options based on SDP and moments relaxation}
\begin{algorithmic}[1]\label{algoSDP}
\REQUIRE Model and payoff parameters, tolerance $\epsilon$
\ENSURE Approximate option price
\STATE   $n=1$, $\mathrm{gap}=1$
\WHILE{\text{$\mathrm{gap} > \epsilon$}}
\STATE Compute $\exp(\tau G_n)$ using Algorithm~\ref{alg:step}
\STATE Compute moments of order $n$ using \eqref{condmoments}
\STATE Solve corresponding $SDP_n$ to get $LowerBound$ and $UpperBound$ 
\STATE  $\mathrm{gap} = |UpperBound - LowerBound|$
\STATE  $n=n+1$
\ENDWHILE
\end{algorithmic}
\end{algorithm}

The following lemma extends the result of Lemma~\ref{lemmanormJ} to the Heston model.
\begin{lemma}
    \label{lem:heston_est}
Let $G_n$ be the matrix representation of the operator $\mathcal{G}$ introduced above with respect to the basis~\eqref{eq:basisJacobi} of $\mathrm{Pol}_{n}(\mathbb{R}^2)$.  
Then the matrix 1-norm of $G_n$ is bounded by
\begin{equation*}
n( r + \kappa + \kappa\theta - \frac{ \sigma^2}{2} ) + \frac12 n^2 ( 1 + |\rho| \frac{\sigma}{2} + \sigma^2).
\end{equation*}
\begin{proof}
Similar to the proof of Lemma \ref{lemmanormJ}.
\end{proof}
\end{lemma}


%% file: chapter4.tex
\section{Numerical experiments}
\label{sec:numerical_experiments}

We have implemented the algorithms described in this paper in {\sc
Matlab} and compare them with Higham's scaling and squaring method
from~\cite{Higham2009}, which typically employs a diagonal Pad\'e
approximation of degree $13$ and is referred to as ``\texttt{expm}''
in the following. The implementation of our algorithms for block
triangular matrices, Algorithm~\ref{alg:full} (fixed scaling
parameter), and Algorithm~\ref{alg:adaptive} (adaptive scaling
parameter), is based on the same scaling and squaring design and are
referred to as ``\texttt{incexpm}' in the following.  All experiments
were run on a standard laptop (Intel Core i5, 2 cores, 256kB/4MB L2/L3
cache) using a single computational thread.

\subsection{Random block triangular matrices}
\label{sec:random}

We first assess run time and accuracy on a randomly generated block
upper triangular matrix $G_n \in \R^{2491\times 2491}$.  There are $46$ diagonal
blocks, of size varying between $20$ and $80$.  The matrix is
generated to have a spectrum contained in the interval $[-80, -0.5]$,
and a well conditioned eigenbasis $X$ ($\kappa_2(X) \approx 100)$.

\begin{figure}[t]
\centering
    \includegraphics[width=0.49\textwidth]{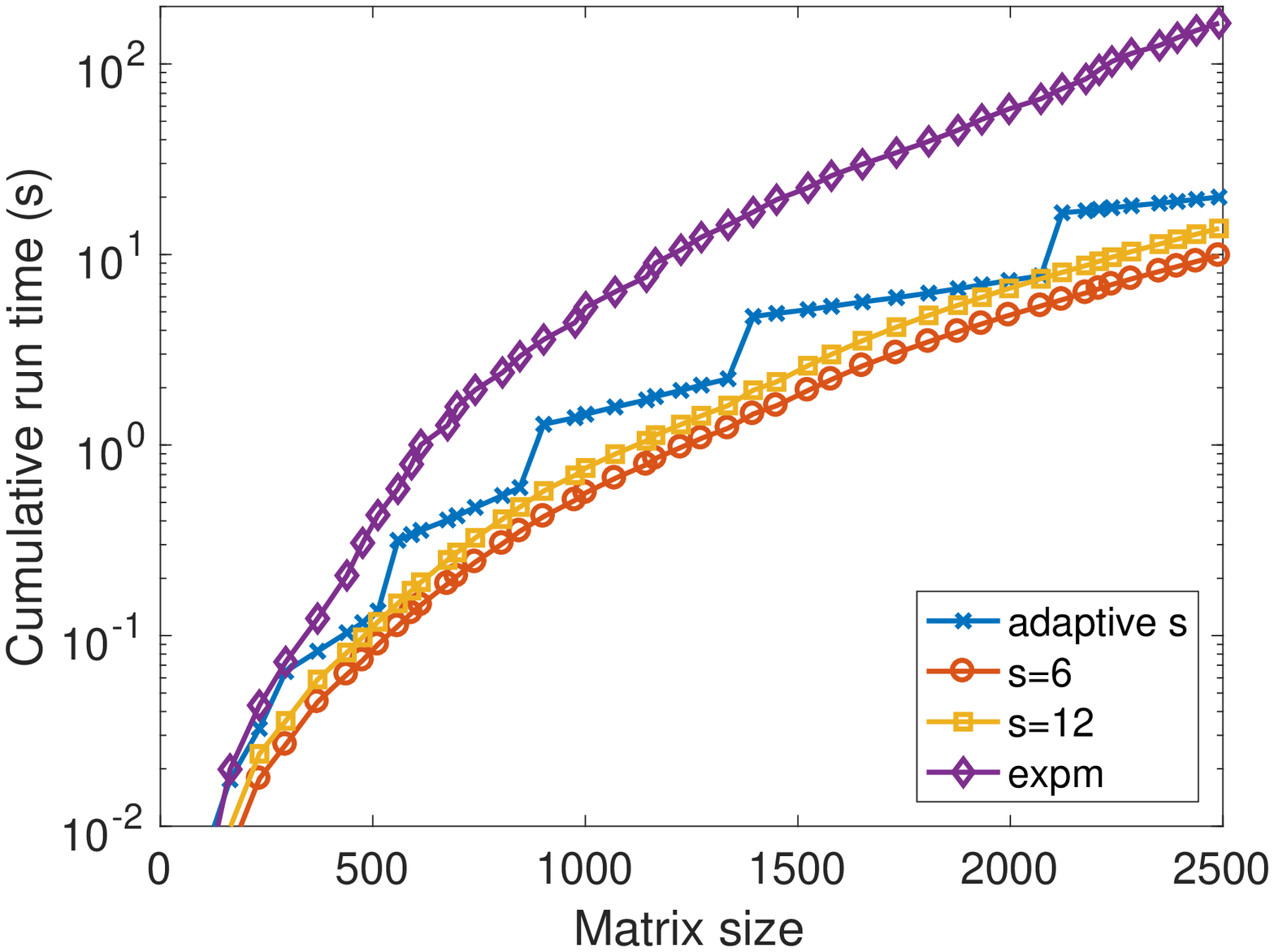}
    \hfill
    \includegraphics[width=0.49\textwidth]{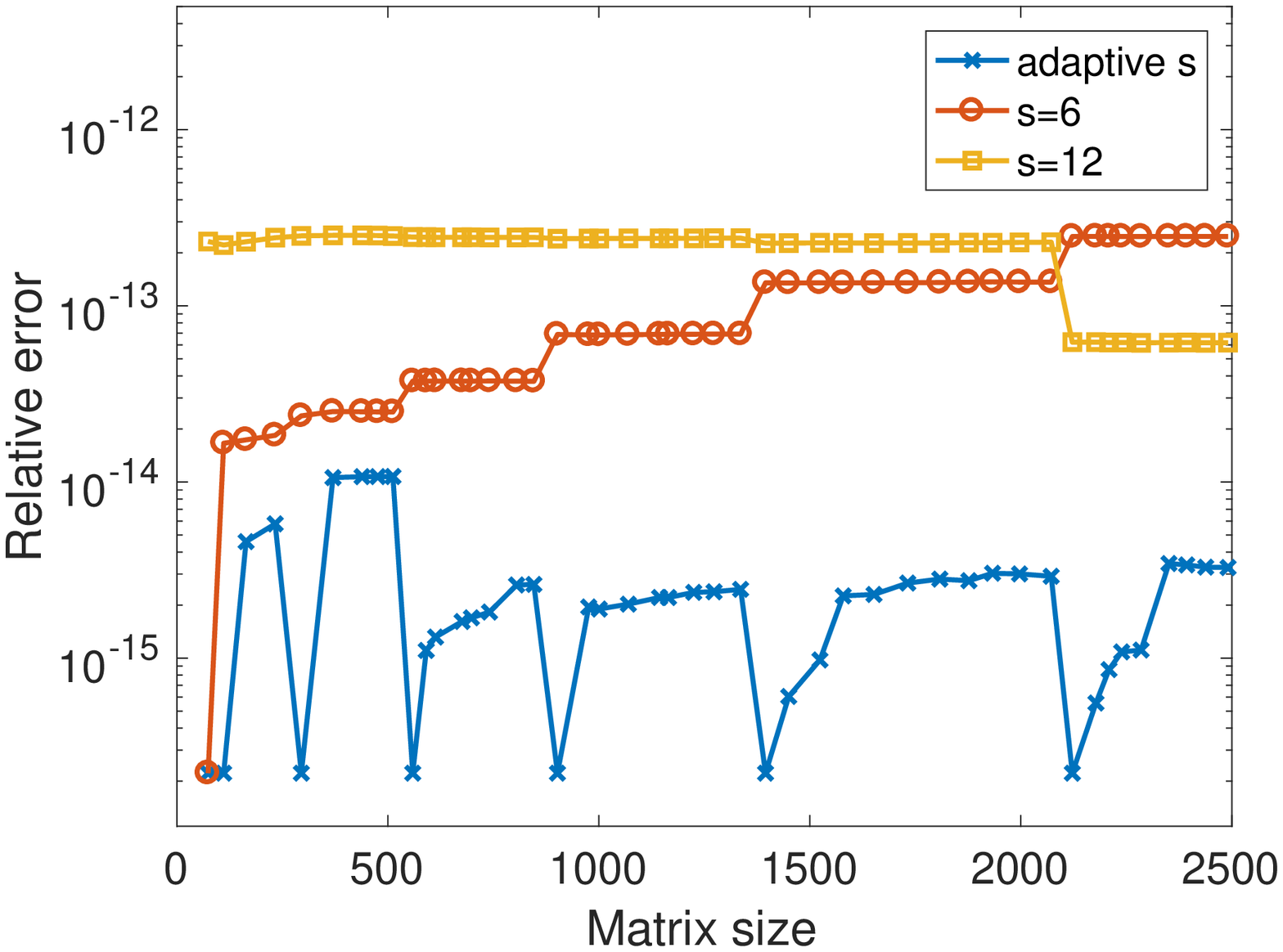}
    \caption{Comparison of \texttt{incexpm} and \texttt{expm} for a
    random block triangular matrix.  \emph{Left:} Cumulative run time
    for computing the leading portions. \emph{Right:} Relative error
    of \texttt{incexpm} w.r.t.~\texttt{expm}.\label{fig:exp_random}}
\end{figure}

Figure~\ref{fig:exp_random} (left) shows the wall clock time for the
incremental computation of all the leading exponentials. 
Specifically, given $0 \le l \le n$, each data point shows the time vs.~$d_l = b_0 + \cdots + b_l$ needed for computing the $l+1$ matrix exponentials $\exp(G_0), \exp(G_1), \dotsc, \exp(G_l)$ when using
\begin{itemize}
 \item \texttt{expm} (by simply applying it to each matrix separately);
 \item \texttt{incexpm} with the adaptive scaling strategy from
Algorithm~\ref{alg:adaptive};
\item \texttt{incexpm} with fixed scaling power $6$ (scaling used by \texttt{expm} for $G_0$);
\item \texttt{incexpm} with fixed scaling power $12$ (scaling used by \texttt{expm} for $G_n$).
\end{itemize}
As expected, \texttt{incexpm} is much faster than naively applying
\texttt{expm} to each matrix separately; the total times for $l = n$
are also displayed in Table~\ref{tab:exp_random}.  For reference we
remark that the run time of \textsc{Matlab}'s \texttt{expm} applied
only the final matrix $G_n$ is $13.65$s, which is very close to the
run time of \texttt{incexpm} with scaling parameter set to $12$ (see
Section~\ref{sec:overall} for a discussion of the asymptotic
complexity).  Indeed, a closer look at the runtime profile of \texttt{incexpm}
reveals that the computational overhead induced by the more
complicated data structures is largely compensated in the squaring
phase by taking advantage of the block triangular matrix structure,
from which \textsc{Matlab}'s \texttt{expm} does not profit
automatically.  It is also interesting to note that the run time of
the adaptive scaling strategy is roughly only twice the run time for
running the algorithm with a fixed scaling parameter $6$, despite its
worse asymptotic complexity.
\begin{table}[t]
\centering
    \caption{Run time and relative error attained by
    \texttt{expm} and \texttt{incexpm} on a random block triangular
    matrix of size $2491$.\label{tab:exp_random}}
    \begin{tabular}{l|r|r}
        Algorithm                  & Time (s) & Rel.~error\\  \hline
        \texttt{expm}              & 163.60   &\\
        \texttt{incexpm} (adaptive)&  20.01   & 3.27e-15\\
        \texttt{incexpm} ($s=6$)   &   9.85   & 2.48e-13\\
        \texttt{incexpm} ($s=12$)  &  13.70   & 6.17e-14\\
    \end{tabular}
\end{table}

The accuracy of the approximations obtained by \texttt{incexpm} is
shown on the right in Figure~\ref{fig:exp_random}.  We assume
\texttt{expm} as a reference, and measure the relative distance
between these two approximations, i.e.,
\begin{equation*}
    \frac{\norm{\texttt{expm}(G_l) - \texttt{incexpm}(G_l)}_F}
        {\norm{\texttt{expm}(G_l)}_F},
\end{equation*}
at each iteration $l$ (quantities smaller than the machine precision
are set to $u$ in Figure~\ref{fig:exp_random}, for plotting purpose).  One notes that the
approximations
of the adaptive strategy remain close to \texttt{expm} throughout the
sequence of computations.  An observed drop of the error down to $u$ for this strategy
corresponds to a restart in Algorithm~\ref{alg:adaptive}; the
approximation at this step is \emph{exactly} the same as the one of
\texttt{expm}.  Even for the fixed scaling parameters $6$ and $12$, the
obtained approximations are quite accurate.

\subsection{Application to option pricing}
\label{sec:exp_jacobi}

We now show results for computing option prices using
Algorithm~\ref{alg:jacobi} for the set of parameters 
\begin{align*}
    &v_0=0.04, \quad  x_0=0,\quad \sigma_w=0.5,\quad \mu_w=0,\quad \kappa=0.5,\quad \theta=0.04,\quad \sigma=0.15, 
    \\ &\rho=-0.5,\quad v_\text{min}=0.01,\quad v_\text{max}=1,\quad r=0,\quad \tau=1/4, \quad k=\log(1.1).
\end{align*}
We use the tolerance $\epsilon = 10^{-3}$ for stopping Algorithm~\ref{alg:jacobi}.

\begin{figure}[t]
    \centering
    \includegraphics[width=0.49\textwidth]{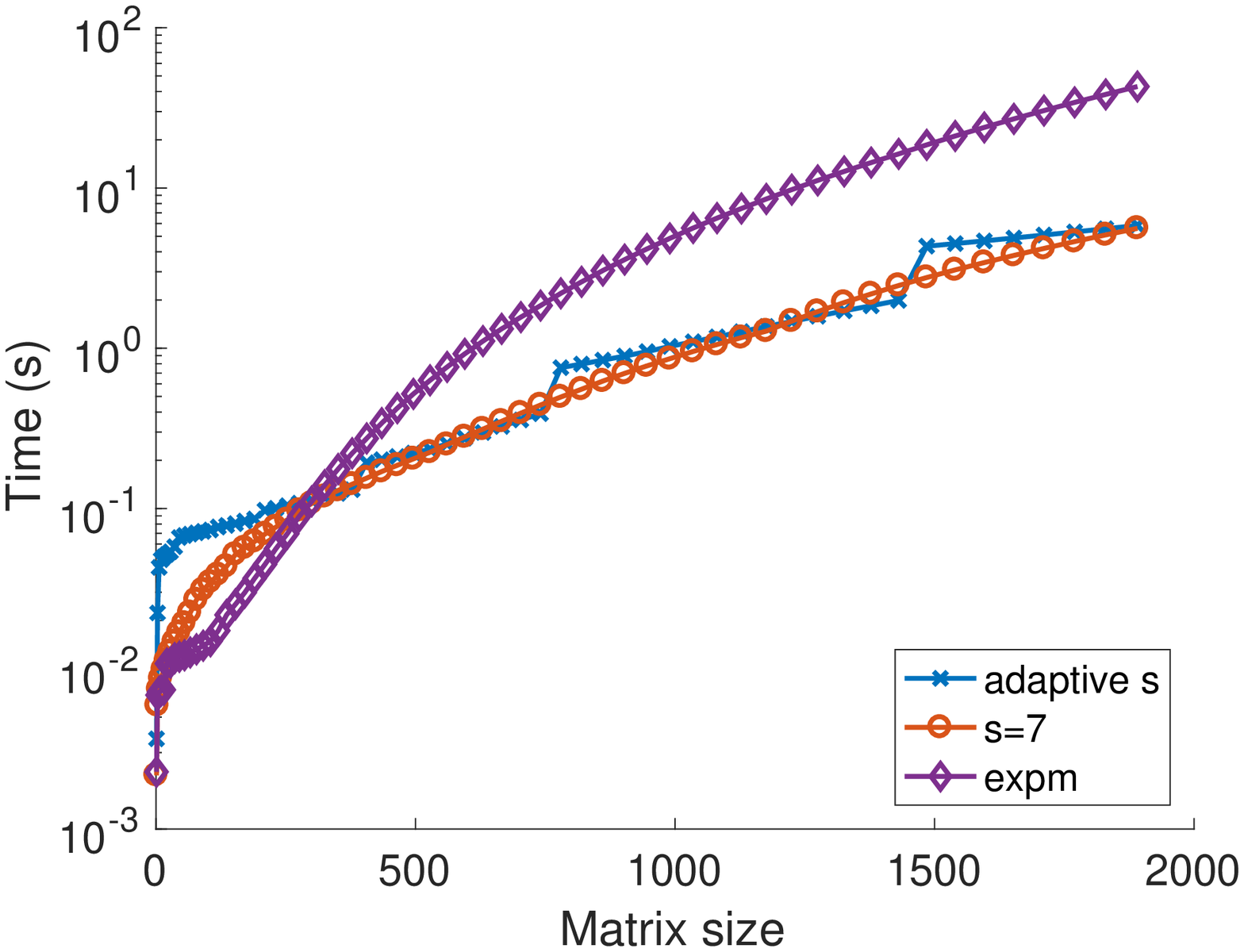}
    \hfill
    \includegraphics[width=0.49\textwidth]{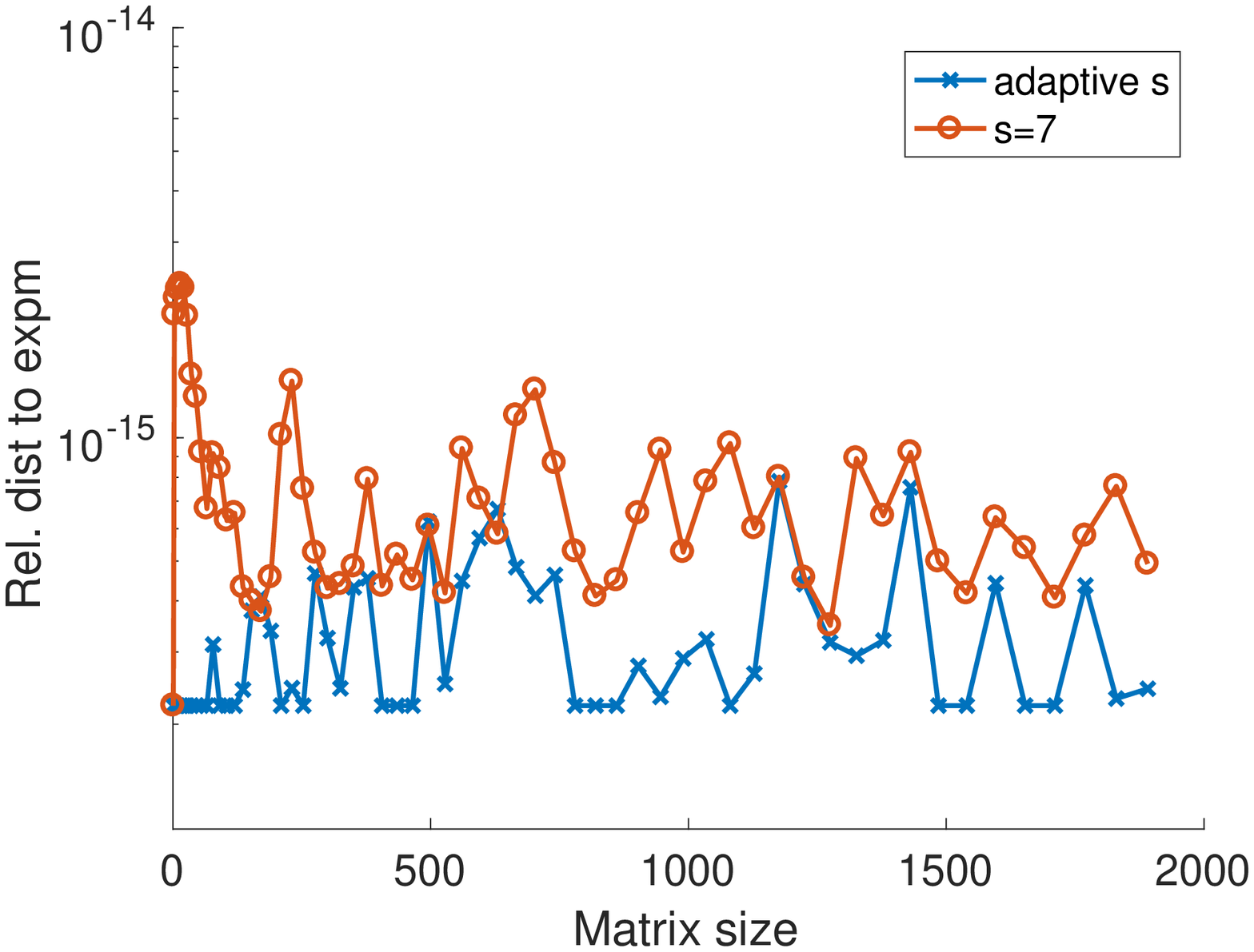}
    \caption{Comparison of \texttt{incexpm} and \texttt{expm} for the block upper triangular matrices arising in the context of the Jacobi model in Algorithm~\ref{alg:jacobi}.
    \emph{Left:} Cumulative run time for computing the leading
    portions. \emph{Right:} Relative error of \texttt{incexpm}
    w.r.t.~\texttt{expm}. \label{fig:exp_jacobi}}
\end{figure}

We explore the use of different algorithms for the computation of the matrix
exponentials in line~\ref{line:jacobiexp} of
Algorithm~\ref{alg:jacobi}: \texttt{incexpm} with adaptive scaling,
\texttt{incexpm} with fixed scaling parameter $s=7$ (corresponding to
the upper bound from Lemma~\ref{lemmanormJ} for $n = 60$), and
\texttt{expm}. Similar to Figure~\ref{fig:exp_random}, the observed
cumulative run times and errors are shown in
Figure~\ref{fig:exp_jacobi}. Again, \texttt{incexpm} is observed to be
significantly faster than \texttt{expm} (except for small matrix
sizes) while delivering the same level of accuracy.  Both
\texttt{incexpm} run times are also close to the run time of
\textsc{Matlab}'s \texttt{expm} applied only to the final matrix $\tau G_n$
($4.64$s).

\begin{table}[t]
    \centering
    \caption{Total run time and option price errors for the
    Jacobi model for $n = 61$. \label{table_jacobi} }
    \begin{tabular}{l|r|r}
        Algorithm & Time (s) & Rel.~price error\\  \hline
        \texttt{expm}               &42.97 & 1.840e-03\\
        \texttt{incexpm} (adaptive) & 5.84 & 1.840e-03\\
        \texttt{incexpm} ($s=7$)    & 5.60 & 1.840e-03\\
    \end{tabular}
\end{table}

Table~\ref{table_jacobi} displays the impact of the different algorithm on the overall Algorithm~\ref{alg:jacobi}, in terms of execution time and accuracy. Concerning accuracy, we computed the relative error with respect to a reference option price computed by considering a truncation order $n = 100$.  It can be observed that there is no difference in accuracy 
for the three algorithms.

\begin{remark}
    \begin{rm}
    The block triangular matrices $G_n$ arising from the generator in
    the Jacobi model actually exhibit additional structure. They are quite sparse and the diagonal blocks
    are in fact permuted triangular matrices (this does not hold for
    polynomial diffusion models in general, though). For example, for $n=2$ the matrix $G_2$ in the Jacobi model is explicitly given by
\begin{equation*}
    G_2 = \left[
    \begin{array}{c|cc|ccc}
        0& r     &\kappa \theta &0&-\frac{\rho \sigma v_\text{max}v_\text{min}}{S}&-\frac{\sigma^2 v_\text{max}v_\text{min}}{S}\\ \hline
        & 0    & 0 &2r&\kappa \theta&0\\
        &  -\frac{1}{2} & -\kappa&1&r+\frac{\rho \sigma (v_\text{max}+v_\text{min})}{S}& 2\kappa\theta +\frac{\sigma^2( v_\text{max}+v_\text{min})}{S}\\ \hline
       &     &&0&0&0\\
        &     & &-1&-\kappa&0\\
        &    & &0&-\frac{1}{2}-\frac{\rho \sigma }{S}&-2\kappa-\frac{\sigma^2}{S}\\
    \end{array}
    \right],
\end{equation*}
    for $S\defby (\sqrt{v_{\max}}-\sqrt{v_{\min}})^2$.
    
    While the
    particular structure of the diagonal blocks is taken into account automatically by \texttt{expm} and
    \texttt{incexpm} when computing the LU decompositions of the diagonal blocks, it is not so easy to benefit from the sparsity.
    Starting from sparse matrix arithmetic, the matrix quickly becomes 
    denser during the evaluation of the initial rational
    approximation, and in particular during the squaring phase.  In
    all our numerical experiments we used a dense matrix
    representation throughout.
    \end{rm}
\end{remark}

We repeated the experiments above for the Heston instead of the Jacobi model, that is, we investigated the impact of using our algorithms for computing the matrix exponentials in Algorithm~\ref{algoSDP}. We found that the results for computing the matrix exponentials themselves look very similar to those for the Jacobi model (Figure~\ref{fig:exp_jacobi}), both
in terms of run time and accuracy, so we refrain from giving further
details here. There is, however, a notable difference. The
evaluation of the stopping criterion requires the solution of two SDPs, which quickly becomes a
computational challenge, eventually completely dominating the time needed for
the computation of the matrix exponentials. 

%% file: chapter5.tex
\section{Summary and future work}

We have presented techniques for scaling and squaring algorithms that
allow for the incremental computation of block triangular matrix
exponentials.  We combined these techniques with an adaptive scaling
strategy that allows for both fast and accurate computation of each
matrix exponential in this sequence (Algorithm~\ref{alg:adaptive}).
For our application in polynomial diffusion models, the run time can
be further reduced by using fixed scaling parameter, determined
through the estimation techniques in Lemmas~\ref{lemmanormJ}
and~\ref{lem:heston_est}.

We observed in our numerical experiments that accurate approximations
to these matrix exponentials can be obtained even for quite small,
fixed scaling parameters.  For the case of two-by-two block triangular
matrices, the results of Dieci and Papini~\cite{Dieci2000,Dieci2001}
support this finding, but an extension of these results to cover a
more general setting would be appreciable.